\documentclass[12pt]{amsart}
\usepackage[left=2.2cm,top=3.3cm,right=2.2cm]{geometry}
\usepackage{amsmath}
\usepackage{amssymb}
\usepackage{setspace}
\usepackage{enumerate}
\begin{document}

\newtheorem*{Lusin}{Lusin's Theorem}

\title{A One-Sentence Inverse Image Proof of Lusin's Theorem}

\begin{abstract}
We give a one-sentence proof of Lusin's theorem. We do not believe our approach, by way of inverse images, is new. However, this particular proof is strikingly clear and succinct.
\end{abstract}

\author[S. Ferguson]{Samuel J. Ferguson}
\address[S. Ferguson]{Courant Institute of Mathematical Sciences, New York University, 251 Mercer St, New York, NY 10012}
\email{ferguson@cims.nyu.edu}

\author[T. Wu]{Tianqi Wu}
\address[T. Wu]{Courant Institute of Mathematical Sciences, New York University, 251 Mercer St, New York, NY 10012}
\email{tw1196@nyu.edu}

\maketitle

Standard proofs of Lusin's theorem using simple functions, cf.\,\cite{bass, rudin}, are sometimes quite elaborate. While they are instructive in learning ``how things are done" in measure theory, we feel that the shorter proof below should be more widely known, and might raise the popularity of this intuitive (cf.\,\cite{royden}) result. The authors do not believe this approach to the proof is new. However, this particular proof appears to be strikingly clear and succinct.

We state Lusin's theorem with the general hypotheses of \cite{rudin} in the Borel measurable context and the conclusion from \cite{bass}. Below, $Y$ is a Borel subset of a locally compact Hausdorff space $X$, and $\mu$ is a measure on the $\sigma$-algebra of Borel subsets of $X$ such that $\mu(K)<\infty$ for each compact subset $K$ of $X$, $\mu(E)=\inf\{\mu(U):E\subset U\text{ and }U\text{ is open}\}$ for each Borel subset $E$, and $\mu(E) = \sup\{\mu(K):K\subset E\text{ and }K\text{ is compact}\}$ for each Borel subset $E$ with $\mu(E)<\infty$.

\begin{Lusin}
Let $Y$ and $\mu$ be as above, with $\mu(Y)<\infty$. If $f:Y\to\mathbb{R}$ is Borel measurable, then for each $\varepsilon > 0$, there exists a subset $C$ of $Y$, which is closed in $Y$'s subspace topology, such that $\mu(Y\setminus C) < \varepsilon$ and the restriction of $f$ to $C$, denoted $f|_{C}$, is continuous.
\end{Lusin}

The idea is that, to ensure $f|_{C}$ is continuous, it suffices to have $f|_{C}^{-1}\left( (a,b)\right)$ open in $C$'s subspace topology for each open interval $(a,b)$ with rational endpoints $a, b$. After all, every open subset of $\mathbb{R}$ is a union of such intervals. As we can use the assumptions on $\mu$ to approximate all of these (countably many) intervals' inverse images at once, we get a candidate for $C$.

\begin{proof}[Proof:]
Given $\varepsilon > 0$, listing all of the open intervals with rational endpoints, $(a_{1}, b_{1})$, \dots, $(a_{n}, b_{n})$, \dots, and using the fact that, by hypothesis, there exist compact $K_{n}$ and open $U_{n}$ such that $K_{n}\subset f^{-1}\left( (a_{n}, b_{n})\right)\subset U_{n}$ and $\mu(U_{n}\setminus K_{n})<\varepsilon/2^{n}$, we can let $C=Y\setminus\cup_{n\geq 1}\left( U_{n}\setminus K_{n}\right)$ and observe that, for each $n\geq 1$,
\[
f|_{C}^{-1}\left( (a_{n}, b_{n})\right) = C\cap f^{-1}\left( (a_{n}, b_{n})\right) = C\cap U_{n},
\]
by definition of $C$.
\end{proof}


\begin{thebibliography}{99}

\bibitem{bass} Bass, R.F. \emph{Real Analysis for Graduate Students}, Second Edition, Createspace, 2013.

\bibitem{royden} Royden, H.L. \emph{Real Analysis}, Third Edition, Macmillan, New York, NY, 1988.

\bibitem{rudin} Rudin, W. \emph{Real and Complex Analysis}, Third Edition, McGraw-Hill, New York, NY, 1987.

\end{thebibliography}
\end{document}